\documentclass[a4paper,10pt]{article}
\usepackage[utf8x]{inputenc}

\usepackage{amssymb}

\newtheorem{theorem}{Theorem}
\newtheorem{lemma}[theorem]{Lemma}
\newtheorem{corollary}[theorem]{Corollary}
\newtheorem{proposition}[theorem]{Proposition}

\def\R{\mathbb{R}}

\def\set#1{\left\{\, #1 \,\right\}}
\def\abs #1{\left| \,#1\, \right|}
\def\norm #1{\left\| \,#1\, \right\|}
\def\qed{\hfill$\square$}

\def\Leg{\mathcal{L}}

\def\calC{\mathcal{C}}

\def\calH{\mathcal{H}}

\def\calX{\mathcal{X}}

\begin{document}

\title{Translation invariance of weak KAM solutions
of the Newtonian n-body problem}

\author{Ezequiel Maderna}

\maketitle
\begin{abstract}
                We consider the Hamilton-Jacobi equation $H(x,d_xu)=c$,
where $c\geq 0$, of the classical $N$-body problem in some Euclidian
space $E$ of dimension at least two. The fixed points of the Lax-Oleinik
semigroup are global viscosity solutions for the critical value of
the constant ($c=0$), also called weak KAM solutions. We show that
all these solutions are invariant under the action by
translations of $E$ in the space of configurations. We also show the
existence of non-invariant solutions for the super-critical equations
($c>0$).
\end{abstract}

\section{Introduction and results}

\subsection{Preliminaries}

Let $E$ be an Euclidean vector space of dimension $k \geq 2$ in which
$N$ punctual masses $m_1,\dots, m_N >0$ are mutually attracted
by the Newtonian gravitational forces. At each time, the position
vectors of the $N$ bodies determine a configuration of the system
$x=(r_1,\dots,r_N)\in E^N$.
The Newtonian potential is the function $U:E^N\to (0,+\infty]$
defined by
\[U(x)=\sum_{1\leq i<j\leq N}\;\frac{m_im_j}{\norm{r_i-r_j}}.\]
It is clear that $U(x)<+\infty$ if and only if $x$ is a
configuration \textsl{without collisions},
meaning that there are no values of $i<j$ such that $r_i=r_j$.
Moreover, if we denote by $\Omega\subset E^N$ the open and
dense subset of configurations without collisions, then the
Newtonian potential is an analytical function on $\Omega$,
and the law of gravitation is equivalent to the second order
differential equation $\ddot x=\nabla U$ on $\Omega$,
where the gradient is taken with respect to the mass
inner product.

As well as the equations of motion, the Hamiltonian and the
Lagrangian of the system admit a synthetic expression in
terms of the mass inner product. We recall that given two
configurations $x=(r_1,\dots,r_N)$ and $y=(s_1,\dots,s_N)$,
their mass inner product is given by
\[x\cdot y = \sum_{i=1}^N\,m_i<r_i,s_i>.\]
With the natural identification $TE^N\simeq E^N\times E^N$,
the Lagrangian function of the system reads
\[L(x,v)=\frac{1}{2}\,\abs{v}^2+U(x)\]
where $\abs{v}=(v\cdot v)^{1/2}$ denotes the
induced norm by the mass inner product on $E^N$.
Consequently, if $p=(p_1,\dots,p_N)$ is in the
dual space $(E^N)^*\simeq (E^*)^N$,
and $\abs{p}$ denotes the operator norm of $p$
for the mass inner product, we have
\[\abs{p}^2=\sum_{i=1}^N\,\frac{1}{m_i}\norm{p_i}^2\]
where $\norm{p_i}$ is the operator norm of $p_i\in E^*$.
Since the Hamiltonian $H$ of the system is defined on
each fiber of the cotangent bundle $T^*E^N$
as the convex dual of the Lagrangian,
we obtain the expression
\[H(x,p)=\frac{1}{2}\,\abs{p}^2-U(x).\]
Recall that here the Legendre transform
$\Leg(x,v)=\partial_vL(x,v)\in T^*E^N$ does
not depends on the configuration $x$, hence
it can be defined as a map
$\Leg:E^N\to (E^N)^*$.
Thus $p=\Leg (v)$ will mean that
$p(w)=w\cdot v$ for all $w\in E^N$.
It is easy to check the Young-Fenchel inequality,
which assures that for all $p\in (E^N)^*$ and
$v\in E^N$ we have $p(v)\leq H(x,p)+L(x,v)$,
and that equality holds if and only if $p=\Leg (v)$.

Finally, we recall that if $M=m_1+\dots+m_N$
is the total mass of the system, the center of mass
of a configuration is given by
the linear map $G:E^N\to E$
\[G(x)=\frac{1}{M}\sum_{i=1}^N\,m_ir_i.\]
For any vector $r\in E$ we have
$x\cdot (r,\dots,r)=<MG(x),r>$,
thus it is clear that the space of configurations
decompose as the orthogonal direct sum
$E^N=\calX \oplus \Delta$, where $\calX=\ker G$
is the subspace of all centered configurations,
and $\Delta$ the diagonal subspace.
We also recall that, since the sum of all forces
acting on the system is zero, the center of mass
of any motion has constant velocity.

\subsection{The Hamilton-Jacobi equation and weak KAM solutions}

Since our mechanical system is autonomous his Hamiltonian
function $H$ is time independent.
Therefore the natural question which
arises is the existence and characterization of the
stationary solutions of the Hamilton-Jacobi equation
$\partial_tS+H(x,\partial_xS)=0$, that is to say,
solutions of the form $S(t,x)=u(x)-ct$.
In this paper we will study the
global stationary solutions,
or in other words, functions $u:E^N\to\R$
satisfying the Hamilton-Jacobi equation $H(x,d_xu)=c$,
more explicitly the equation
\begin{equation}
\label{HJ}
\frac{1}{2}\sum_{i=1}^N
\frac{1}{m_i}\norm{\frac{\partial u}{\partial r_i}}^2-
\sum_{i<j}\frac{m_im_j}{\norm{r_i-r_j}}=c.
\end{equation}
Considering the mutual distances between the bodies
can be arbitrarily large, it is clear that equation
(\ref{HJ}) only has global solutions for values of
the constant $c\geq 0$.

A classical global solution, or
global smooth solution,
is a continuous function $u:E^N\to\R$
whose restriction to the open and dense set
of configurations without collisions $\Omega$
is of class $C^1$, and such that (\ref{HJ})
is verified at each point of this set.
Until now, the only known smooth global solutions are
solutions for the case $N=2$, that is to say,
for the Kepler problem in the space $E$. A direct
computation shows that functions
$u(r_1,r_2)=\pm\norm{r_1-r_2}^{1/2}$
are global smooth solutions.

Recently, the author has shown in \cite{Mad}
the existence
of weak KAM solutions for a family of
$N$-body problems with homogeneous potentials,
including the Newtonian case.
They are global viscosity solutions of the
critical equation $H(x,d_xu)=0$.
The existence was obtained by application
of a fixed point theorem to the action
of the Lax-Oleinik semigroup on the convex
set of viscosity subsolutions.

Moreover, if $\Gamma$ is a given symmetry group
of the system, then the linear space of all
invariant functions has nontrivial intersection
with the set of subsolutions
since clearly constant functions are
invariant subsolutions.
Considering that the semigroup preserves
the set of all invariant subsolutions,
which is also convex,
the existence of invariant weak KAM solutions
is also obtained (theorem 3 in \cite{Mad}).
It is evident that the application of a given
isometry of the Euclidian space $E$ to each component
of a configuration defines a symmetry of the system.
Thus the translation group of $E$, as well as
the orthogonal group $O(E)$, acts naturally
on the space of configurations by symmetries.
Depending on the values of the masses,
other symmetry groups can act on the space
of configurations. For example,
in the non-generic case in which the masses
of some bodies are equal, every permutation
of the position vectors of these bodies
are symmetries.

We shall remark that for a smooth, convex
and superlinear Hamiltonian function $H:T^*M\to\R$
on a compact manifold $M$, weak KAM
solutions are always invariant under the
action of the identity component of a compact
symmetry group (see \cite{Maderna2002}).

Even if the the
orthogonal group is compact, there are
always non-invariant solutions for the
action of his identity component $SO(E)$.
The critical Busemann functions explicitly
computed in \cite{Mad} for the Kepler problem
in the plane with equal masses, are examples
of such non-invariant weak KAM solutions.
Nevertheless, weak KAM solutions are always
invariant under translations.

\begin{theorem}\label{thm}
Every weak KAM solution of the Newtonian $N$-body
problem in an Euclidean space $E$ of dimension
at least two is translation invariant.
More precisely, given a weak KAM
solution $u:E^N\to \R$, we have that
$$u(r_1,\dots,r_N)=u(r_1+r,\dots,r_N+r)$$ for every configuration
$x=(r_1,\dots,r_N)\in E^N$ and every $r\in E$.
\end{theorem}

The hypothesis $\dim E\geq 2$ will be only used
in order to apply Marchal's theorem
(whose proof in the case $\dim E>3$ is due to
D.~Ferrario and S.~Terracini) to the calibrated
curves of a given solution.
This said, there is no evidence that suggest
the theorem is not valid in the one dimensional case.

We will show that the invariance under translations of
a weak KAM solution is equivalent to the fact that
all his calibrated curves are motions with fixed center
of mass (proposition \ref{prop}). On the other
hand, a very strong property of the calibrated curves
of a weak KAM solution is the free time minimization
property, which we recall below.
Therefore, the proof of theorem \ref{thm} will
follow from lemma \ref{ftm}, where Marchal's
theorem is applied. The dynamics
when $t\to +\infty$ of the free time minimizers
of the Newtonian $N$-body problem is described
in \cite{DaLuzMad}, where is proved
that such kind of motions are completely parabolic
and asymptotic to some special central configurations.

The following corollary will be easily deduced from
the theorem. It shows as a counterpart, that
the super-critical equations (\ref{HJ}), that is
for the values of $c>0$, can be non-invariant under
translations.

\begin{corollary}\label{coro}
If $u_0:E^N\to\R$ is a weak KAM solution
of the $N$-body problem, and $r\in E$ is an arbitrary
vector, then the function $$u(x)=u_0(x)+<G(x),r>$$
is a global viscosity solution of the
Hamilton-Jacobi equation
$$H(x,d_xu)=\frac{1}{2}\norm{r}^2.$$
\end{corollary}

The following section is devoted to prove the above results.
It is interesting to remark that all the statements and
proofs remains valid if we replace the Newtonian potential
by homogeneous potential of degree $\alpha\in (-2,0)$.

\section{Proofs}

We start recalling some basic notions of the calculus
of variations that will be used. The Lagrangian action
of an absolutely continuous curve
$\gamma:[a,b]\to E^N$
\[A(\gamma)=\int_a^b\,
L(\gamma(t),\dot\gamma(t))\,dt\]
is a well defined function
of $\gamma$ taking values in $(0,+\infty]$.
For two given configurations $x,y\in E^N$,
the set of all absolutely continuous curves
starting at $x$ and arriving to $y$ in time
$T>0$ will be denoted $\calC (x,y,T)$.
It is well known that the Lagrangian action
is lower semicontinuous, which implies that
the infimum of its restriction to any set
$\calC (x,y,T)$ is always reached. In that
follows $\phi(x,y,T)$ will denote the minimum
of the Lagrangian action restricted
to $\calC (x,y,T)$.

Let us recall Marchal's theorem on
the fixed endpoint problem (cf. \cite{Chen},\cite{FerTer}).
It says that minimizers of the Lagrangian
action in a set $\calC (x,y,T)$ not
suffers collisions in any interior time.
In other words, if an absolutely continuous
curve $\gamma:[a,b]\to E^N$ is such that
$A(\gamma)\leq A(\sigma)$ for every
$\sigma\in\calC (\gamma(a),\gamma(b),b-a)$,
then the curve $\gamma$ avoid collisions
except perhaps in the prescribed endpoints,
that is to say,
$\gamma(t)\in\Omega$ for all $t\in (a,b)$.
In particular, the restriction of a such
curve to the open interval $(a,b)$ satisfies
the equation of motion, hence it must be
of class $C^\infty$. Furthermore, if
one of its endpoints, say $\gamma(a)$,
is a configuration without collisions,
then the solution can be extended to an
interval of the form $(a-\epsilon,b)$
for some $\epsilon>0$,
which implies that $\dot\gamma(a)$ exists
and coincide with
$\lim_{t\to a^+}\dot\gamma(t)$.

\subsection{Action potential, free time minimizers
and calibrated curves}

We will recall some results from \cite{Mad} as
well as some well known features of weak KAM
theory.

\subsubsection{The action potential}
The set of all absolutely continuous curves starting
at $x$ and ending at $y$ without restriction
of time will be denoted $\calC (x,y)$.
The action potential,
also called the Ma\~{n}\'{e}
critical potential in Aubry-Mather theory,
is the function $\phi:E^N\times E^N\to \R^+$
\[\phi(x,y)=
\inf\set{A(\gamma)\mid \gamma\in\calC (x,y)}
=\inf_{T>0}\phi(x,y,T).\]
For convenience, we will use the
maximum norm in $E^N$,
\[\norm{(r_1,\dots,r_N)}=
\max_{i=1,\dots ,N} \norm{r_i}_E.\]
The action potential is a distance
function in the space of configurations.
Moreover, it is H\"{o}lder continuous
in each variable
with respect to any norm in $E^N$.
In fact, there is a constant $\eta>0$
such that
\[\phi(x,y)\leq \eta \norm{x-y}^{1/2}\]
for every pair of configurations. The constant
$\eta$ only depends on the number of bodies $N$
and the total mass $M$.
Since several of the here considered notions
are easily characterized in terms of the
action potential, it will play a central role.

Another important property of the action potential
that we will also use, is that it is locally Lipschitz
in $\Omega\times\Omega$.
More precisely, if $C\subset\Omega$ is a compact set,
then there is a constant $K>0$ such that
$\phi(x,y)\leq K\,\norm{x-y}$ whenever
$(x,y)\in C\times C$. Of course the Lipschitz constant
grows as well as the compact set $C$ approaches
some given collision configuration.

\subsubsection{Free time minimizers}
We call a free time minimizer, or semistatic curve,
an absolutely continuous curve
such that the action of each compact
segment equals the action potential between
the endpoints of the segment. In other words,
a curve $\gamma:I\to E^N$ is a free time minimizer
if we have
$A(\gamma\mid_{[a,b]})=\phi(\gamma(a),\gamma(b))$
whenever $[a,b]\subset I$ is a compact subinterval.
When $\dim E\geq 2$, Marchal's theorem imlpies that
a free time minimizer defined on an interval $[a,b]$
avoid collisions for $t\in(a,b)$.

\subsubsection{Viscosity subsolutions as dominated functions}
On the other hand we will be concerned with
the notion of weak solution of (\ref{HJ})
introduced by M.~Crandall and P.-L.~Lions
in \cite{CranLions},
namely the notion of viscosity solution.
The classical way to define them is
in terms of test functions, as
functions which are at the same time
viscosity subsolutions and viscosity
supersolutions.
It is well known that the set
of viscosity subsolutions of the critical
Hamilton-Jacobi equation $H(x,d_xu)=0$,
which we call the set of dominated functions,
is exactly the convex set
\[\calH=\set{u:E^N\to\R \mid
u(y)-u(x)\leq \phi(x,y)\textrm{ for all } x,y\in E^N}.\]
By the previous considerations, we know
in particular that
functions in $\calH$ are H\"{o}lder continuous
and locally Lipschitz in $\Omega$. Since $\Omega$
is a set of total measure, by
Rademacher's theorem we also know that they
are differentiable almost everywhere.

\subsubsection{Calibrated curves}
For a given dominated function $u\in\calH$, we say that
an absolutely continuous curve $\gamma:I\to E^N$
is calibrated for $u$, if for each compact subinterval
$[a,b]\subset I$ we have that
$u(\gamma(b))-u(\gamma(a))=A(\gamma\mid_{[a,b]})$.
Note that a calibrated curve is always a free time minimizer.
On the other hand, if $u\in\calH$ and $\gamma$ is a 
calibrated curve for $u$, defined on $I=(a,b)$
and of class $C^1$, then $u$ is differentiable
at $\gamma(t)$ for all $t\in (a,b)$, and
$d_{\gamma(t)}u=\Leg(\dot\gamma(t))$.
It follows, using the Young-Fenchel inequality,
that $H(\gamma(t),d_{\gamma(t)}u)=0$ for all
$t\in (a,b)$.

\subsubsection{Characterization of weak KAM solutions}
The action of the Lax-Oleinik semigroup $(T_t)_{t\geq 0}$
is well defined in $[0,+\infty)\times\calH$ by
\[T_tu(x)=\inf\set{u(y)+\phi(y,x,t)\mid y\in E^N}\]
for $t>0$ and $T_0=id_\calH$.
A weak KAM solution is a fixed point of the Lax-Oleinik semigroup,
that is to say, a function $u\in\calH$ such that
$T_tu=u$ for all $t\geq 0$. 
It is proved in \cite{Mad} that such fixed points exist and
that in fact, they are global viscosity solutions of the
critical Hamilton-Jacobi equation.
Moreover they are characterized between dominated functions
as follows: $u\in\calH$ is a weak KAM solution if and only
if for each configuration $x\in E^N$ there is a calibrated
curve $\gamma_x:[0,+\infty)\to E^N$ such that $\gamma_x(0)=x$.
By all the above considerations, we deduce that if $u$ is
a weak KAM solution
differentiable at a configuration $x\in\Omega$, then
there is only one calibrated curve, since it must coincide
with the motion generated by the initial conditions
$\gamma(0)=x$ and $\dot\gamma(0)=\Leg^{-1}(d_xu)$.

\subsection{Proof of theorem \ref{thm}}

As we have said, the proof is a direct consequence of the following
lemma and proposition.

\begin{lemma}\label{ftm}
If $\dim E\geq 2$ and $x:[0,+\infty)\to E^N$ is a free time minimizer
then the center of mass $G(x(t))$ is constant.
\end{lemma}

\textsl{Proof.}
First, we must observe that Marchal's theorem
implies $x(t)\in\Omega$
for all $t>0$.
Therefore $x$ satisfies the equation of motion
on $(0,+\infty)$, from which we deduce that $G(x(t)$
has constant velocity, even at $t=0$.
With the notation $x_0=x(0)$ and $G_0=G(x_0)$, we
have that $G(x(t))$ is of the form $G_0+tv$
for some vector $v\in E$.
We will prove that $v=0$.

Let $y:[0,+\infty)\to E^N$ be the internal motion
associated with $x$, that is to say, the curve
defined by $x(t)=y(t)+\delta(G_0+tv)$, where
$\delta:E\to E^N$ is the map $\delta(r)=(r,\dots,r)$.
It is clear that for all $t\geq 0$ we
have that $U(x(t))=U(y(t))$ and $G(y(t))=0$.
On the other hand, a computation shows that
$\dot x\cdot \dot x=\dot y\cdot \dot y +M\norm{v}^2$
for all $t>0$.  
Hence we deduce that for all $T>0$ we have
\[A\left(x\mid_{[0,T]}\right)=
A\left(y\mid_{[0,T]}\right)+
\frac{1}{2}\,T\,M\norm{v}^2.\]
Since $x$ is a free time minimizer, and $\phi$ verify
the triangle inequality, we also obtain
\[A\left(x\mid_{[0,T]}\right)=\phi(x_0,x(T))\leq
\phi\left(x_0,x_0+T\delta(v)\right)+
\phi\left(x_0+T\delta(v),x(T)\right).\]
Let $y_T$ be the curve obtained
translating the restriction of
$y$ to the interval $[0,T]$
in such a way that the
center of mass of $y_T$ becomes fixed at $G(x(T))$.
Thus, $y_T$ is the curve defined for $t\in [0,T]$
by $y_T(t)=y(t)+\delta(G_0+Tv)$. Both curves
have the same action, but now we have
$\gamma_T\in\calC(x_0+T\delta(v),x(T))$.
Therefore,
\[\phi\left(x(0)+T\delta(v),x(T)\right)\leq
A\left(y\mid_{[0,T]}\right).\]
We use now the H\"{o}lder estimation for
the action potential, and we get the bound
\[\phi\left(x(0),x(0)+T\delta(v)\right)\leq
\eta\norm{T\delta(v)}^{1/2}=
T^{1/2}\;\eta\norm{v}^{1/2}.\]
We can conclude that the inequality
\[\frac{1}{2}\,T\,M\norm{v}^2\;\leq
\;T^{1/2}\;\eta\norm{v}^{1/2}\]
holds for all $T>0$, but this is possible
only in case $v=0$.
\qed

\begin{proposition}\label{prop}
A weak KAM solution $u:E^N\to\R$ is invariant under
translations if and only if all calibrated curves have
constant center of mass.
\end{proposition}

\textsl{Proof.}
We start showing the necessity of the condition.
Let $u$ be a weak KAM solution which is invariant under
translations. At each point of differentiability
$x\in E^N$ we must have have $\Delta\subset \ker d_xu$.
Let $x\in E^N$ be any configuration,
and let $\gamma_x:[0,+\infty)$ be a calibrated curve for
$u$ starting at $x$. We know that $u$ is differentiable
at $\gamma_x(t)$ for all $t>0$ and that 
$d_{\gamma_x(t)}=\Leg (\dot\gamma_x(t))$,
that is to say, $d_xu(w)=w\cdot \dot\gamma_x(t)$
for all $w\in E^N$. It follows that
$\dot\gamma_x(t)\in \Delta^\perp=\calX$
for $t>0$, hence that
$G(\gamma_x(t))$ is constant.

Let now $u$ be a weak KAM solution such that
all calibrated curves have constant center of mass.
Let $x\in E^N$ be some given configuration, and
choose a non zero vector $r\in E$.
We will prove that $u(x)=u(x+\delta(r))$, where
$\delta(r)=(r,\dots,r)\in \Delta$.

Since $u$ is continuous and the set $\Omega$ is dense
and translation invariant, it suffices to consider
configurations $x\in\Omega$.
Let $$B(\epsilon,r)=
\set{z\in E^N\mid z\cdot \delta(r)=0,\;
\norm{z}\leq\epsilon}$$ be the closed ball
of radius $\epsilon>0$
in the orthogonal complement
of $\delta(r)\in\Delta$. Since
$x\in \Omega$, we can assume that $\epsilon$
is so small that the set
$$B=\set{x+z\mid z\in B(\epsilon,r)}$$ is
contained in $\Omega$, which implies
that also the compact cylinder
\[C=\set{y+t\delta(r)\mid
y\in B,t\in[0,1]}\] is contained in $\Omega$.
Let us call $C'\subset C$ the set
of points where $u$ is differentiable.
Recall that $u$ is a dominated function,
hence it is a Lipschitz function on $C$,
and $C'$ has total measure in $C$.
For $q\in C'$, let $\gamma_q:[0,+\infty)\to E^N$
be the only calibrated curve for $u$
with $\gamma_q(0)=q$. We know that
$d_qu(w)=w\cdot \dot\gamma_q(0)$
for all $w\in E^N$.
By using the hypothesis we deduce that
$\dot\gamma_q(0)\in\calX=\ker G=\Delta^\perp$,
and therefore that $\Delta\subset \ker d_qu$
for all $q\in C'$.

We consider now the Lipschitz continuous function
$f:B\times [0,1]\to \R$ defined by
$f(y,t)=u(y+t\delta(r))-u(y)$. It is clear that
$\partial_tf$ is zero almost everywhere.
Therefore, applying Fubini's theorem $\partial_tf$
on a cylider of the form $A\times [0,1]$, where
$A\subset B$ is an open set, we obtain
\[\int_A f(y,1)\,dy=0.\]
Since the open subset $A\subset B$
is arbitrary, we must have $f(y,1)=0$ for all $y\in B$.
Hence, we have $f(x,1)=0$,
that is to say, $u(x+\delta(r))=u(x)$.
\qed

\subsection{Viscosity solutions of
the supercritical equations}

\textsl{Proof of corollary \ref{coro}.}
Let $u_0$ be a weak KAM solution and $r\in E$.
We will use the fact that $u_0$ is a global
solution of $H(x,d_xu)=0$ in the
viscosity sense.
As before, we call $\delta:E\to \Delta$
the trivial isomorphism. Thus function
$u_r$ can be written as $u_r=u_0+a_r$, where
$a_r$ is the linear function
$a_r(x)=x\cdot \delta(r)$.

Let $x_0$ be a given configuration, and
suppose that $\varphi\in C^1(E^N)$ is
such that $\varphi-u_r$ has a local minimum
at $x_0$.
Since $u_0$ is a viscosity solution of
$H(x,d_xu)=0$
and we have that $(\varphi-a_r)-u_0$
has a local minimum at $x_0$, we have
$H(x_0,d_{x_0}(\varphi-a_r))\leq 0$.

On the other hand,
by theorem \ref{thm} we know
that $u_0$ is translation invariant, hence
$\Delta\subset \ker d_{x_0}(\varphi-a_r)$.
Since
$d_{x_0}a_r=a_r$, it is clear that
$\calX\subset \ker d_{x_0}a_r$ and
that $\abs{d_{x_0}a_r}^2=M^*\norm{r}^2$.
Therefore we deduce that
\[\abs{d_{x_0}\varphi}^2=
\abs{d_{x_0}(\varphi-a_r)}^2+
M^*\norm{r}^2\]
and we conclude that
\[H(x_0,d_{x_0}\varphi)=
H(x_0,d_{x_0}(\varphi-a_r))+
\frac{1}{2}M^*\norm{r}^2\leq
\frac{1}{2}M^*\norm{r}^2,\]
proving that $u_r$ is a viscosity subsolution
of the supercritical equation at $x_0$.
With a symmetric argument it can be proved
that $u_r$ is a viscosity supersolution of
the same equation.
\qed

\bibliographystyle{amsplain}

\begin{thebibliography}{00}

\bibitem{DaLuzMad}
  A. Da Luz and E. Maderna,
  \textsl{On the free time minimizers
  of the Newtonian $N$-body problem},
  preprint, http://premat.fing.edu.uy/2009.htm


\bibitem{Chen}
  A. Chenciner,
  \textsl{Action minimizing solutions of the
  Newtonian $n$-body problem:
  from homology to symmetry},
  Proceedings of the ICM, Vol.~III
  (Beijing, 2002),  279--294,
  Higher Ed. Press, Beijing, 2002.

\bibitem{CranLions}
  M.G. Crandall and P.-L. Lions,
  \textsl{Viscosity Solutions of Hamilton-Jacobi Equations},
  Trans. Amer. Math. Soc.
  \textbf{277} (1983), 1--42.
  

\bibitem{FerTer}
  D. Ferrario and S. Terracini,
  \textsl{On the existence of collisionless equivariant
  minimizers for the classical $n$-body problem},
  Invent. Math. \textbf{155} (2004), no. 2, 305--362.

\bibitem{Maderna2002}
  E. Maderna,
  \textsl{Invariance of global solutions of the
  Hamilton-Jacobi equation}, Bull. Soc. Math. France
  \textbf{130} (2002) no.4, 493--506.

\bibitem{Mad}
  E. Maderna,
  \textsl{On weak KAM theory for $N$-body problems},
  Ergodic Theory Dynam. Systems.
  Available on CJO 2011 doi:10.1017/S0143385711000046

\end{thebibliography}

\end{document}